\documentclass[12pt]{article}
\usepackage{amsmath,amssymb}
\usepackage{fullpage}

\newcommand{\R}{\mathbb R}
\newcommand{\C}{\mathbb C}

\newcommand{\Z}{\mathbb Z}
\newcommand{\T}{\mathbb T}

\renewcommand{\O}{\Omega}
\renewcommand{\liminf}{\mathop{\hbox{\rm lim inf}}\limits}

\newcommand{\supp}{{\rm supp}}

\newcommand{\A}{{\cal A}}
\newcommand{\cL}{{\cal L}}
\newcommand{\cC}{{\cal C}}
\newcommand{\cB}{{\cal B}}
\newcommand{\cT}{{\cal T}}
\newcommand{\cI}{{\cal I}}
\newcommand{\cP}{{\cal P}}
\newcommand{\cR}{{\cal R}}
\newcommand{\Linfg}{L^\infty (G)}

\newcommand{\spec}{{\rm spec}}
\newcommand{\cM}{{\cal M}}
\newcommand{\cMT}{{\cal M}_T(\Sigma)}

\newtheorem{defin}{Definition}[section]
\newtheorem{thm}[defin]{Theorem}
\newtheorem{ex}[defin]{Example}
\newtheorem{prop}[defin]{Proposition}
\newtheorem{lemma}[defin]{Lemma}

\newtheorem{rem}[defin]{Remarks}

\begin{document}

\title{Transference in Spaces of Measures}

\author{Nakhl\'e H.\ Asmar\footnote{
Department of Mathematics,
University of
Missouri,
Columbia, Missouri 65211  U.\ S.\ A.
}, 
 Stephen J.\ Montgomery--Smith\footnote{
Department of Mathematics,
University of
Missouri,
Columbia, Missouri 65211  U.\ S.\ A.
}, and Sadahiro Saeki\footnote{
Department of Mathematics,
Kansas State University,
Manhattan, Kansas 66506  U.\ S.\ A.
}
\date{}
}
\maketitle
			     %
%%%%%%%%%%%%%%%%%%%%%%%%%%%%%%%%%%%%%%%%%%%%%
%%%%%%%%%%%%%%%%%%%%%%%%%%%%%%%%%%
\section{Introduction}

\renewcommand{\thefootnote}{\roman{footnote}}

%%%%%%%%%%%%%%%%%%%%%%%%%%%%%%%%%%%%%%%%%%%%%
%%%%%%%%%%%%%%%%%%%%%%%%%%%%%%%%%%
%%%%%%%%%%%%%%%%%%%%%%%%%%%%%%%%%%%%%%%%%%%%%
%%%%%%%%%%%%%%%%%%%%%%%%%%%%%%%%%%

Transference theory for $L^p$ spaces is a powerful tool with
many fruitful applications to singular integrals, ergodic
theory, and spectral theory of operators \cite{cal,cw}.
These methods afford a unified approach to many problems
in diverse areas, which before were proved by a variety of methods.

The purpose of this paper is to bring about a similar approach to spaces of measures.  Our main transference result is 
motivated by the extensions of the classical F.\&M.~Riesz Theorem 
due to Bochner
 \cite{bs},
Helson-Lowdenslager \cite{hl1, hl2}, 
de Leeuw-Glicksberg \cite{deleeuwglicksberg}, 
Forelli \cite{forelli}, and others.  It might seem
that these extensions should all be obtainable 
via transference methods,
and indeed, as we will show, these are exemplary 
illustrations
of the scope of our main result.

It is not straightforward to extend the classical transference methods
of Calder\'on, Coifman and Weiss to spaces of measures.  
First, their methods make use of 
averaging techniques and the amenability of the group of representations.
The averaging techniques simply do not work with
measures, and do not preserve analyticity.  Secondly, and most importantly,
their techniques require that the representation is strongly continuous.
For spaces of measures, this last requirement would be prohibitive, 
even for the simplest representations such as translations.  Instead, we
will introduce a much weaker requirement, which we will call `sup path
attaining'.  By working with sup path attaining representations, we are
 able to prove a new transference principle with interesting applications.
For example, we will show how to derive with ease generalizations of Bochner's theorem and Forelli's main result.  The Helson-Lowdenslager theory, concerning representations of groups with ordered dual groups, 
is also within reach, but it will
be treated in a separate paper.

Throughout $G$ will denote a 
locally compact abelian group with dual group
$\Gamma$.  
The symbols $\Z$, $\R$ and $\C$ 
denote the integers, the real and complex numbers, respectively.  
If $A$ is a set, we denote the indicator
function of $A$ by $1_A$.
For $1\leq p<\infty$,
the space of Haar measurable functions $f$ on $G$ with
$\int_G|f|^p dx<\infty$ will be denoted by 
$L^p(G)$.  The space of essentially
bounded functions on $G$ will be denoted by 
$L^\infty(G)$.  The expressions ``locally null''
and ``locally almost everywhere'' will have the same meanings as
in \cite[Definition (11.26)]{hr1}.

Our measure theory is borrowed from
\cite{hr1}.
In particular, the space of all complex regular Borel measures 
on $G$,
denoted by $M(G)$, consists of all complex measures  
arising from bounded linear functionals on $\cC_0(G)$,
the Banach space of continuous functions on
$G$ vanishing at infinity. 
 
Let $(\O, \Sigma)$ denote a
measurable space, 
where $\O$ is a set and $\Sigma$ is a
sigma algebra of subsets of $\O$.  Let $M(\Sigma)$ denote the 
Banach space of complex measures on $\Sigma$ with the
total variation norm, and let
$\cL^\infty(\Sigma)$ denote the space of measurable
bounded functions on $\Omega$.  

Let
$T:\ t\mapsto T_t$ denote a representation of $G$
by isomorphisms of $M(\Sigma)$.  
We suppose that $T$ is uniformly bounded,
i.e., there is a positive constant $c$ such that
for all $t\in G$, we have
%%%%%%%%%%%
%%%%%%%%%%%
%%%%%%%%%%
\begin{equation}
\|T_t\|\leq c .
\label{uniformlybded}
\end{equation}
%%%%%%%%%%%%%%%%%%%%%%%%%%%%%%%%%%%%%%%%%%%%%%%%%%%%%%%%%
%%%%%%%%%%%%%%%%%%%%%%%%%%%%%%%%%%%%%%%%%%%%%%%%%%%%%%%%
%%%%%%%%%%%%%%%%%%%%%%%%%%%%%%%%%%%%%%%%%%%%%%%%%%%%%%%%%
\begin{defin}
A measure $\mu\in M(\Sigma)$ is 
called weakly measurable (in symbols, $\mu\in{\cal M}_T(\Sigma)$)
if for every $A\in \Sigma$ 
the mapping $t\mapsto T_t\mu(A)$ is 
Borel measurable on $G$.  
\label{weakmble}
\end{defin}
%%%%%%%%%%%%%
%%%%%%%%%%%%%%%%%%%%%%%%%%%%%%%%%%%%%%%%%%%%%%%%%%%%%%%%%
%%%%%%%%%%%%%%%%%%%%%%%%%%%%%%%%%%%%%%%%%%%%%%%%%%%%%%%%
%%%%%%%%%%%%%%%%%%%%%%%%%%%%%%%%%%%%%%%%%%%%%%%%%%%%%%%%
%%%%%%%%%%%%%%%%%%%%%%%%%%%%%%%%%%%%%%%%%%%%%%%%%%%%%%%%%
%%%%%%%%%%%%%%%%%%%%%%%%%%%%%%%%%%%%%%%%%%%%%%%%%%%%%%%%
%%%%%%%%%%%%%%%%%%%%%%%%%%%%%%%%%%%%%%%%%%%%%%%%%%%%%%%%%
%%%%%%%%%%%%%%%%%%%%%%%%%%%%%%%%%%%%%%%%%%%%%%%%%%%%%%%%%
%%%%%%%%%%%%%%%%%%%%%%%%%%%%%%%%%%%%%%%%%%%%%%%%%%%%%%%%

Given a measure $\mu\in \cMT$ and a Borel measure 
$\nu \in M(G)$, we define the `convolution' 
$\nu*_T\mu$ on $\Sigma$ by
\begin{equation}
\nu*_T\mu (A)=\int_G T_{-t}\mu(A) d\nu(t)
\label{Tconv}
\end{equation}
for all $A\in\Sigma$.

We will assume throughout this paper
that the representation $T$ commutes with the 
convolution (\ref{Tconv}) in the following sense:
 for each $t\in G$, 
%%%%%%%%%%%%%%%%%%%%%%%%%%%%%%%%%%%%%%%%%%%%%%%%%%%%%%%%%
%%%%%%%%%%%%%%%%%%%%%%%%%%%%%%%%%%%%%%%%%%%%%%%%%%%%%%%%%
\begin{equation} 
 T_t(\nu*_T\mu)=\nu*_T(T_t\mu).
\label{commut}
\end{equation}
%%%%%%%%%%%%%%%%%%%%%%%%%%%%%%%%%%%%%%%%%%%%%%%%%%%%%%%%%
%%%%%%%%%%%%%%%%%%%%%%%%%%%%%%%%%%%%%%%%%%%%%%%%%%%%%%%%%
%%%%%%%%%%%%%%%%%%%%%%%%%%%%%%%%%%%%%%%%%%%%%%%%%%%%%%%%%
%%%%%%%%%%%%%%%%%%%%%%%%%%%%%%%%%%%%%%%%%%%%%%%%%%%%%%%%%%%%
%%%%%%%%%%%%%%%%%%%%%%%%%%%%%%%%%%%%%%%%%%%%%%%%%%%%%%
%%%%%%%%%%%%%%%%%%%%%%%%%%%%%%%%%%%%%%%%%%%%%%%%%%%%%%%%%%%%%%
%%%%%%%%%%%%%%%%%%%%%%%%%%%%%%%%%%%%%%%%%%%%%%%%%%%%
%%%%%%%%%%%%%%%%%%%%%%%%%%%%%%%%%%%%%%%%%%%%%%%%%%%%%%%%%
Condition (\ref{commut}) holds if, for example,
for all $t\in G$, the adjoint of $T_t$ maps 
$\cL^\infty(\Sigma)$ into itself.  In symbols,
%%%%%%%%%%%
%%%%%%%%%%%
%%%%%%%%%%
\begin{equation}
T_t^*: \cL^\infty(\Sigma) \rightarrow \cL^\infty(\Sigma).
\label{adjointT}
\end{equation}
%%%%%%%%%%%%%%%%%%%%%%%%%%%%%%%%%%%%%%%%%%%%%%%%%%%%%%%%%
%%%%%%%%%%%%%%%%%%%%%%%%%%%%%%%%%%%%%%%%%%%%%%%%%%%%%%%%%
%%%%%%%%%%%%%%%%%%%%%%%%%%%%%%%%%%%%%%%%%%%%%%%%%%%%%%%%%
%%%%%%%%%%%%%%%%%%%%%%%%%%%%%%%%%%%%%%%%%%%%%%%%%%%%%%%%%
For proofs we refer the reader to 
\cite{ams2}.  Using (\ref{uniformlybded}) and 
(\ref{commut}), it can be shown that $\nu*_T\mu$ is a measure in $\cMT$,
%%%%%%%%%%%%%
%%%%%%%%%%%%%%%%%%%%%%%%%%%%%%%%%%%%%%%%%%%%%%%%%%%%%%%%%
%%%%%%%%%%%%%%%%%%%%%%%%%%%%%%%%%%%%%%%%%%%%%%%%%%%%%%%%
\begin{equation}
\|\nu*_T\mu\|\leq  c\|\nu\|\|\mu\|,
\label{normofconv}
\end{equation}
%%%%%%%%%%%%%
%%%%%%%%%%%%%%%%%%%%%%%%%%%%%%%%%%%%%%%%%%%%%%%%%%%%%%%%%
where $c$ is as in 
(\ref{uniformlybded}), and 
%%%%%%%%%%%%%%%%%%%%%%%%%%%%%%%%%%%%%%%%%%%%%%%%%%%%%%%%%%%%%%%%%%%%%
%%%%%%%%%%%%%%%%%%%%%%%%%%%%%%%%%%%%%%%%%%%%%%%%%%%%%%%%%
%%%%%%%%%%%%%%%%%%%%%%%%%%%%%%%%%%%%%%%%%%%%%%%%%%%%%%%%
%%%%%%%%%%%%%%%%%%%%%%%%%%%%%%%%%%%%%%%%%%%%%%%%%%%%%%%%%
%%%%%%%%%%%
%%%%%%%%%%%
%%%%%%%%%%
\begin{equation} 
\sigma*_T(\nu*_T\mu)=(\sigma*\nu)*_T\mu,
\label{associative}
\end{equation}
%%%%%%%%%%%%%%%%%%%%%%%%%%%%%%%%%%%%%%%%%%%%%%%%%%%%%%%%%
%%%%%%%%%%%%%%%%%%%%%%%%%%%%%%%%%%%%%%%%%%%%%%%%%%%%%%%%%
%%%%%%%%%%%%%%%%%%%%%%%%%%%%%%%%%%%%%%%%%%%%%%%%%%%%%%%%%
%%%%%%%%%%%%%%%%%%%%%%%%%%%%%%%%%%%%%%%%%%%%%%%%%%%%%%%%%%%%%%%%%%%
%%%%%%%%%%%%%%%%%%%%%%%%%%%%%%%%%%%%%%%%%%%%%%%
%%%%%%%%%%%%%%%%%%%%%%%%%%%%%%%%%%%%%%%%%%%%%%%%%%%%%%%%%%%%%%%%%%%
%%%%%%%%%%%%%%%%%%%%%%%%%%%%%%%%%%%%%%%%%%%%%%%
%%%%%%%%%%%%%%%%%%%%%%%%%%%%%%%%%%%%%%%%%%%%%%%%%%%%%%%%%
for all $\sigma , \nu \in M(G)$ and $\mu \in \cMT$
(see \cite{ams2}).

We come now to a definition which is fundamental to our work.
%%%%%%%%%%%
%%%%%%%%%%%
%%%%%%%%%%%%%%%%%%%%%%%%%%%%%%%%%%%%%%%%%
%%%%%%%%%%%%%%%%%%%%%%%%%%%%%%%%%%%%%%%%%%%%%%%%%
\begin{defin}
A representation $T=(T_t)_{t\in G}$ of a locally compact abelian group
$G$ in $M(\Sigma)$ is said to be sup path attaining
if it is uniformly bounded, satisfies property (\ref{commut}), and if there 
is a constant $C$ such that for every weakly 
measurable $\mu\in {\cal M}_T(\Sigma)$ we have
\begin{equation}
\| \mu\| \leq C\sup \left\{
{\rm ess\ sup}_{t\in G} \left|
\int_\O h d (T_t\mu)
\right|
:\ \ h\in \cL^\infty(\Sigma),\ \|h\|_\infty\leq 1
\right\}.
\label{ineqhypa}
\end{equation}
\label{hypa}
\end{defin}

The fact that the mapping $t\mapsto \int_\O h d (T_t\mu)$ is 
measurable is a simple consequence of
the measurability of the mapping $t\mapsto T_t\mu(A)$ for every
$A\in\Sigma$.  

Examples of  sup path attaining representations 
 will be presented in the following section.  
Proceeding toward the main result of this paper, we 
recall some basic definitions from spectral theory.

%%%%%%%%%%%%%%%%%%%%%%%%%%%%%%%%%%%%%%
%%%%%%%%%%%%%%%%%%%%%%%%%%%%%%%%%%%%%%%%%%%%%%%%%%%%%%%%%%%
 
If $I$ is an ideal in $L^1(G)$, let
$$Z(I)=\bigcap_{f\in I}
\left\{
\chi\in\Gamma:\ \ \widehat{f}(\chi)=0
\right\}.$$
The set $Z(I)$ is called the zero set of $I$.
For a weakly measurable $\mu\in M(\Sigma)$, let
$$\cI (\mu)=\{f\in L^1(G):\ \ f*_T\mu =0\}.$$
Using properties of the convolution $*_T$,
it is straightforward to show that $\cI(\mu)$ is a closed ideal
in $L^1(G)$.
%%%%%%%%%%
%%%%%%%%%%%%
%%%%%%%%%%%%%%
%%%%%%%%%%%%%%%%
\begin{defin}
The $T$-spectrum of a 
weakly measurable $\mu\in \cMT$ is defined by 
%%%%%%%%%%
%%%%%%%%%%
%%%%%%%%%%
\begin{equation}
{\rm spec}_T (\mu)= \bigcap_{f\in \cI(\mu)}
\left\{
\chi\in\Gamma:\ \ \widehat{f}(\chi)=0
\right\}=Z(\cI(\mu)).
\label{specsbt}
\end{equation}
%%%%%%%%%%
%%%%%%%%%%
%%%%%%%%%%
%%%%%%%%%%
\label{Tspectrum}
\end{defin}
%%%%%%%%%%
%%%%%%%%%%%%
%%%%%%%%%%%%%%
If $S\subset \Gamma$, let
$$L_S^1=L_S^1(G)=\left\{f\in L^1(G):\ \widehat{f}=0\ \mbox{outside of}\ S\right\}\,.$$
Our transference result concerns convolution operators on
$L^1_S(G)$ where $S$ satisfies
a special property, described as follows.
\begin{defin}
A subset $S\subset\Gamma$ is a $\cT$-set if, given any compact 
$K\subset S$, each neighborhood of $0\in\Gamma$ contains a nonempty open set
$W$ such that $W+K\subset S$. 
\label{t-set}
\end{defin}

\begin{ex}
{\rm (a)  If $\Gamma$ is a locally compact abelian group, then any open subset of $\Gamma$ is a $\cT$-set.  In particular, if $\Gamma$ is discrete then every subset of $\Gamma$ is a $\cT$-set.\\
(b)  The set $\left[ a,\infty\right. )$ is a $\cT$-subset of $\R$, for all $a\in\R$.  \\
(c)  Let $a\in\R$ and $\psi:\ \Gamma \rightarrow \R$ be a continuous homomorphism.  Then 
$S=\psi^{-1}([a,\infty))$ is a $\cT$-set.\\
(d)  Let $\Gamma=\R^2$ and $S=\{(x,y):\ y^2\leq x\}$.  
Then $S$ is a $\cT$-subset of $\R^2$ such that there is no nonempty open set $W\subset \R^2$ such that $W+S\subset S$. 
}
\end{ex}

The main result of this paper is the following transference
theorem.
%%%%%%%%%%%%%%%%%%%%%%%%%%%%%%%%%%%%%%%%%%%%%%%%%%%%%%
%%%%%%%%%%%%%%%%%%%%%%%%%%%%%%%%
%%%%%%%%%%%%%%%%%%%%%%%%%%%%%%%%%%%%%%%%%%%%%%%%%%%
%%%%%%%%%%%%%%%%%%%%%%%%%%%
\begin{thm}
Let $T$ be a sup path attaining representation of a locally compact abelian group
$G$ by isomorphisms of $M(\Sigma)$ and let $S$ be a $\cT$-subset of
$\Gamma$.
Suppose that $\nu$ is a measure in $M(G)$ such that
%%%%%%%%%%%
%%%%%%%%%%%%
\begin{equation}
\|\nu*f\|_1\leq \|f\|_1
\label{hyptransference2}
\end{equation}
%%%%%%%%%%%
%%%%%%%%%%%%
for all $f$ in $L_S^1(G)$.  Then for every 
weakly measurable $\mu \in M(\Sigma)$ with
${\rm spec}_T ( \mu )\subset S$
we have
%%%%%%%%%%%
%%%%%%%%%%%%
\begin{equation}
\|\nu*_T\mu\|\leq c^3 C \|\mu\|,
\label{conctransference2}
\end{equation}
%%%%%%%%%%%%%%%%%%%%%%%%%%%%%%%%%%%%%%%%%%%
%%%%%%%%%%%%%%%%%%%%%%%%%%%%%%%%%%%%%%%%%%%%
%%%%%%%%%%%%%%%%%%%%%%%%%%%%%%%%%%%%%%%%%%%%
%%%%%%%%%%%%%%%%%%%%%%%%%%%%%
where $c$ is as in (\ref{uniformlybded}) and 
$C$ is as in (\ref{ineqhypa}).
\label{maintransferenceth}
\end{thm}

To state Forelli's main result in \cite{forelli}, let us recall 
two definitions of Baire sets.  
Suppose that $\Omega$ is
a topological space. 
Usually,
the collection of Baire sets,\ $\cB_0=\cB_0(\O)$, is defined 
as the sigma algebra generated by sets that are compact and also countable
intersections of open sets.  A second definition is to define $\cB_0$
 as the minimal
sigma algebra so that compactly supported continuous functions are measurable.
For locally compact Hausdorff
topological spaces, these two
definitions are equivalent.

  Suppose that $(T_t)_{t\in \R}$ is a group of 
homeomorphisms
of $\Omega$ onto itself such that the mapping $(t,\omega) \mapsto T_t 
\omega$
is jointly continuous.  The maps $T_t$ induce isomorphisms $T_t$ on
the space of Baire measures via the identity $T_t \mu(A) = \mu(T_t(A))$.

If $\nu$ is a Baire measure, we will say that it is 
quasi-invariant if $T_t\nu$ and $\nu$ are mutually absolutely
continuous for all $t \in \R$, that is, for all $A \in \cB_0$ we have
that $|\nu|(A) = 0$ if and only if $|T_t\nu|(A) = 0$.  If $\mu$ is
a Baire measure, we will say that $\mu$ is $T$-analytic if
$$ \int_{\R} T_t\mu(A) h(t) dt = 0, $$
for all $A \in \cB_0$ and all $h \in H^1(\R)$, where
$$H^1(\R)=\left\{f\in L^1(\R):\ \widehat{f}(s)=0\ \mbox{for all}\ s\leq 0\right\}\,.$$ 
The main result of
Forelli \cite{forelli} says the following.

\begin{thm}
\label{forelli}
Let $\Omega$ be a locally compact Hausdorff 
topological space, and let $(T_t)_{t\in \R}$ 
be a group of homeomorphisms of $\Omega$ onto itself such that the maps 
$(t,\omega) \mapsto T_t \omega$ is jointly continuous.  Suppose that
$\mu$ is a $T$-analytic Baire measure, and that $\nu$ is a quasi-invariant
Baire measure.  Then both $\mu_s$ and $\mu_a$ are $T$-analytic, where $\mu_s$
and $\mu_a$ are the singular and absolutely continuous parts of $\mu$ with respect
to $\nu$.
\end{thm}

The proofs of both this result in \cite{forelli}, and also its predecessor
\cite{deleeuwglicksberg},
are long and difficult.  Furthermore, it is hard to understand
why one must use Baire measures instead of Borel measures.
As it turns out, the mystery of why we need to use Baire measures in Forelli's work is reduced to the fact that
such representations as described above on the Baire measures are sup path attaining.  By working with sup path attaining representations,
we are able to prove a more general version of Theorem \ref{forelli}.
We do not need $\Omega$ to be locally compact Hausdorff.  
Or
we might suppose that $(T_t)_{t\in \R}$ is any group of uniformly bounded
isomorphisms satisfying (\ref{commut}) on any Lebesgue space (that is, a 
countably generated sigma algebra).  

%%%%%%%%%%%%%%%%%%%%%%%%%%%%%%%%%%%%%%%%%%%%%%%%%%%%%%%%%%%%%%%%%%%%%%%%%%%%%%%%%%%%%%%%%%
%%%%%%%%%%%%%%%%%%%%%%%%%%%%%%%%%%%%%%%%%%%%%%%%%%%%%%%%%%%%%%%%%%%%%%%%%%%%%%%%%%%%%%%%%%
%%%%%%%%%%%%%%%%%%%%%%%%%%%%%%%%%%%%%%%%%%%%%%%%%%%%%%%%%%%%%%%%%%%%%%%%%%%%%%%%%%%%%%%%%%%
%%%%%%%%%%%%%%%%%%%%%%%%%%%%%%%%%%%%%%%%%%%%%
%%%%%%%%%%%%%%%%%%%%%%%%%%%%%%%%%%%%%%%%%%%%%%%%%
%%%%%%%%%%%%%%%%%%%%%%%%%%%%%%%%%%%%%%%%%%%%%%%%%%%%%%%%%%%%%%%%%%%%%%%%%%%%%%%%%%%%%%%%%%%%%%%%%
%%%%%%%%%%%%%%%%%%%%%%%%%%%%%%%%%%%%%%%%%%%%%%%%%%%%%%%%%%%%%%%%%%%%%%%%%%%%%%%%%%%%%%%%%%%%%%%%%
%%%%%%%%%%%%%%%%%%%%%%%%%%%%%%%%%%%%%%%%%%%%%%%%%%%%%%%%%%%%%%%%%%%%%%%%%%%%%%%%%%%%%%%%%%%%%%%%%
%%%%%%%%%%%%%%%%%%%%%%%%%%%%%%%%%%%%%%%%%%%%%%%%%%%%%%%%%%%%%%%%%%%%%%%%%%%%%%%%%%%%%%%%%%%%%%%%%%
%%%%%%%%%%%%%%%%%%%%%%%%%%%%%%%%%%%%%%%%%%%%%%%%%%%%
%%%%%%%%%%%%%%%%%%%%%%%%%%%%%%%%%%%%%%%%%%%%%%%%%
%%%%%%%%%%%%%%%%%%%%%%%%%%%%%%%%%%%%%%%%%%%%%%%%%%%%%%%%%%%%%%%%%%%%%%%%%%%%%%%%%%%%%%%%%%%%%%%%%%
%%%%%%%%%%%%%%%%%%%%%%%%%%%%%%%%%%%%%%%%%%%%%%%%%%%%%%%%%%%%%%%%%%%%%%%%%%%%%%%%%%%%%%%%%%%%%%%%%%%
%%%%%%%%%%%%%%%%%%%%%%%%%
%%%%%%%%%%%%%%%%%%%%%%%%%%%%%%%%%%%%%%%%%%%%%%%%%%%%%%%%%%%%%%%%%%%%%%%%%%%%%%%%%%%%%%%%%%%%%%%%%%%
%%%%%%%%%%%%%%%%%%%%%%%%%%%%%%%%%%%%%%%%%%%%%%%%%%%%%%%%%%%%%%%%%%%%%%%%%%%%%%
%%%%%%%%%%%%%%%%%%%%%%%%%%%%%%%%%%%%%%%%%
\section{Sup Path Attaining Representations}

%%%%%%%%%%%%%
%%%%%%%%%%%%
%%%%%%%%%%%

We first note that for a sup path attaining representation $T$ of $G$ we
have
\begin{equation}
\| \mu\| \leq C\sup \left\{
{\rm ess\ sup}_{t\in G} \left|
\int_\O h d (T_t\mu)
\right|
:\ \ h {\rm \ is\ a\ simple\ function\ with\ }\ \|h\|_\infty\leq 1
\right\}.
\label{ineqhypa'}
\end{equation}
To see this, note that for any bounded measurable function $h$,
there exist a sequence of simple functions $s_n$ that converge uniformly
to $h$.  Then it is easy to see that
$$  {\rm ess\ sup}_{t\in G} \left|
\int_\O s_n d (T_t\mu)
\right|
\to
{\rm ess\ sup}_{t\in G} \left|
\int_\O h d (T_t\mu)
\right| ,$$
and from this (\ref{ineqhypa'}) follows easily.

\

Our first example
is related to the setting of Forelli \cite{forelli}.
%%%%%%%%%%%%%
%%%%%%%%%%%%
%%%%%%%%%%%%%
%%%%%%%%%%%%
%%%%%%%%%%%
\begin{ex}
{\rm  Let $G$ be a locally compact abelian group.
Suppose that $\Omega$ is a topological
space and $\left( T_t\right)_{t\in G}$ is a
 group of homeomorphisms of $\Omega$ onto itself
such that the mapping
$$(t,\omega)\mapsto T_t\omega$$
is jointly continuous.
Suppose that $\A$ is an algebra of bounded continuous complex valued
functions
on $\Omega$ such that if $h \in \A$,
and if $\varphi:\C\to\C$ is any bounded continuous function,
and if $t \in G$,
then $\varphi\circ h\circ T_t\in \A$.  Let $\sigma(\A)$ denote the
minimal sigma-algebra so that
functions from $\A$ are measurable.
For any measure $\mu\in M(\sigma(\A))$, and $A\in \sigma(\A)$,
define
$T_t\mu(A)=\mu(T_t(A))$, where $T_t(A)=\{T_t\omega:\ \omega\in A\}$.
Note that $T$ satisfies
(\ref{uniformlybded}) and (\ref{commut}).
To discuss the weak measurability of $\mu$, and the sup path attaining
property of $T$, we need that for each $h \in \A$ that the map
$t \mapsto \int_\Omega h \, d(T_t\mu)$ is continuous.  
This 
crucial property follows
for any measure $\mu\in M(\sigma(\A))$, by the dominated convergence theorem, if, for example, $G$ is metrizable.  For arbitrary locally compact abelian groups, it is enough to require that $\mu$ has sigma-compact support.

Henceforth, we assume that $G$ is metrizable, or that $\mu\in M(\sigma(\A))$ has sigma-compact support, and proceed to show that
$\mu$ is
weakly measurable in the sense of Definition \ref{weakmble}, and that
the representation $T$ is sup path attaining.
Let
%%%%%%%%%%%%%%%%%%%%%%%%%%%%%%%%%%%%%%%%%%%%%
%%%%%%%%%%%%%%%%%%%%%%%%%%%%%%%%%%%%%%%%%%%%%
$$\cR=\{A\in \sigma(\A) :\ 1_A=\lim_n f_n,\ f_n \in \A,\
\|f_n\|_\infty\leq 1\},$$
%%%%%%%%%%%%%%%%%%%%%%%%%%%%%%%%%%%%%%%%%%%%%
%%%%%%%%%%%%%%%%%%%%%%%%%%%%%%%%%%%%%%%%%%%%%
and let
$$\cC=\{A\in \sigma(\A):  \ t\mapsto T_t\mu(A)\
{\rm is\ Borel\ measurable}\}.$$
Clearly, $\cC$ is a monotone class, closed under
nested unions and intersections.  Also, $\cR$
is an algebra of sets, closed under finite unions and
set complementation.  Furthermore, it is clear that $\cR \subset \cC$.
Hence, by the monotone class theorem,
it follows that $\cC$ contains the
sigma algebra generated by $\cR$.
Now, if $f\in\A$ and $V \subset \C$ is open, then $f^{-1}(V) \in \cR$, because
$$1_{f^{-1}(V)}(\omega)=
\lim_{n\rightarrow \infty} \min \{ n\ {\rm dist}(f(\omega),\C\setminus V)
,1\}.$$
Consequently, we have that $\cC = \sigma(\A)$, that is,
$\mu$ is weakly measurable.

Next, let us show that $\A$ is dense in $L^1(|\mu|)$.  Let
$g:\Omega\to\C$ be a bounded $\sigma(\A)$-measurable
function such that $\int_\Omega hg\,d\mu = 0$ for all
$ h \in \A$.  Define $\cR$ as above, and let
$$ \cC = \left\{ A \in \sigma(\A): \int_A g \, d\mu = 0 \right\} .$$
Then $\cC$ is a monotone class containing $\cR$, and so arguing as
before, $\cC = \sigma(\A)$, that is, $g = 0$ almost everywhere with respect
to $|\mu|$.  Thus it follows by the Hahn-Banach Theorem that
$\A$ is dense in $L^1(|\mu|)$.  Hence
$$ \|\mu\|
   =
   \sup\left\{\left|\int_\Omega h \, d\mu \right|: h \in \A,\ \|h\|_\infty \le 1\right\}
   \leq
   \sup\left\{{\rm ess\ sup}_{t\in G}
       \left| \int_\Omega h \, d(T_t\mu)\right| : h \in \A,\ \|h\|_\infty \le 1\right\} ,$$
where the last inequality follows from the fact that 
the map $t \mapsto \int_\Omega h \, d(T_t\mu)$ is continuous.
Hence, $T$ is sup path attaining with $C = 1$.
}
\label{exhypa1}
\end{ex}
%%%%%%%%%%%%%%%%%%%%%%%%%%%%%%%%%%%%
%%%%%%%%%%%%%%%%%%%%%%%%%%%%%%%%%%%%
Taking $\A$ to be the uniformly continuous functions on a group,
we then derive the following useful example.
%%%%%%%%%%%%%%%%%%%%%%%%%%%%%%%%%%%%
%%%%%%%%%%%%%%%%%%%%%%%%%%%%%%%%%%%%
\begin{ex}
{\rm
Suppose that $G_1$ and $G_2$ are locally compact abelian groups and that
 $\phi:\ G_2\rightarrow G_1$ is a continuous homomorphism.
Define an action of $G_2$ on $M(G_1)$ (the regular Borel measures
on $G_1$) by translation by $\phi$.
Hence, for $x\in G_2, \mu\in M(G_1)$, and any
Borel subset $A\subset G_1$, let $T_x\mu(A)=\mu(A+\phi(x))$.
Then every $\mu\in M(G_1)$ is weakly measurable, the
representation
is sup path attaining with constants $c = 1$ and $C = 1$.
}
\label{exhypa2}
\end{ex}
%%%%%%%%%%%%%%%%%%%%%%%%%%%%%%%%%%%%
%%%%%%%%%%%%%%%%%%%%%%%%%%%%%%%%%%%%
%%%%%%%%%%%%%%%%%%%%%%%%%%%%%%%%%%%%
In the following example no topology is required on
the measure space.
%%%%%%%%%%%%%%%%%%%%%%%%%%%%%%%%%%%%
%%%%%%%%%%%%%%%%%%%%%%%%%%%%%%%%%%%%
\begin{ex}
{\rm
Let $(X,\Sigma)$ be an abstract Lebesgue 
space, that is, $\Sigma$ is countably generated.  
Then 
any uniformly bounded representation $T$ of $G$ 
by isomorphisms of $M(\Sigma)$ is sup path attaining.  
To see this, note that since $\Sigma$ is countably 
generated, there is a countable subset ${\cal A}$ 
of the unit ball of ${\cal L}^\infty(\Sigma)$ such that for any 
$\mu\in M(\Sigma)$ 
we have 
$$
\|\mu\| =\sup\left\{
\left| \int_X h d \mu\right| :\ \ h\in 
{\cal A}
\right\} .
$$
If $\mu$ is weakly measurable, then
for $g\in {\cal A}$ and for locally almost all $u\in G$, we have
\begin{equation}
\left|
\int_X g d (T_u \mu)
\right|
\leq {\rm ess\ sup}_{t\in G}
\left|
\int_X
g d (T_t \mu)
\right|.
\label{want1}
\end{equation}
Since ${\cal A}$ is countable we can find a subset $B$ 
of $G$ such that the complement of $B$ is locally null, 
and such that (\ref{want1}) holds for every $u\in B$ 
and all $g\in {\cal A}$.  For $u\in B$, 
take the sup in (\ref{want1}) over all $g\in {\cal A}$, and get
$$\| T_u \mu\| \leq 
\sup_{g\in {\cal A}}
{\rm ess\ sup}_{t\in G}
\left|
\int_X
g d (T_t \mu)
\right|.
$$
But since $\|\mu\|\leq c \| T_u \mu\|$, it follows that 
$T$ is sup path attaining with $C=c$.}
\label{exhypa3}
\end{ex}
%%%%%%%%%%%%%
%%%%%%%%%%%%
%%%%%%%%%%%
%%%%%%%%%%%%%
%%%%%%%%%%%%
%%%%%%%%%%%
%%%%%%%%%%%%
%%%%%%%%%%%

Sup path attaining representations satisfy the
following property, which was introduced in \cite{ams2} and was
called hypothesis $(A)$. 
%%%%%%%%%%%%%
%%%%%%%%%%%%
%%%%%%%%%%%
\begin{prop}
Suppose that $T$ is sup path attaining 
and $\mu$ is weakly measurable such that for 
every $A\in \Sigma$ we have
$$T_t\mu(A)=0$$
for locally almost all $t\in G$.  Then $\mu=0$.
\label{weakhypa}
\end{prop}
%%%%%%%%%%%%%
%%%%%%%%%%%%
%%%%%%%%%%%
The proof is immediate and follows from 
(\ref{ineqhypa'}).\\

The key in Proposition \ref{weakhypa}
is that the set of $t\in G$ 
for which $T_t\mu(A)=0$ depends on $A$. 
If this set were the same for all
$A\in \Sigma$, then the conclusion of the
proposition would trivially hold for
any representation by isomorphisms of $M(\Sigma)$. 

For further motivation, we recall the following 
example from \cite{ams2}.
%%%%%%%%%%%%%
%%%%%%%%%%%%
%%%%%%%%%%%
\begin{ex}
{\rm
(a) Let $\Sigma$ denote the sigma algebra of 
countable and co-countable subsets of $\R$.  
Define $\nu\in M(\Sigma)$ by
$$
\nu (A)=\left\{
\begin{array}{ll}
1 & \mbox{if $A$ is co-countable,}\\
0 &  \mbox{if $A$ is countable.}
\end{array}
\right.
$$
Let $\delta_t$ denote the point mass at $t\in \R$, 
and take $\mu=\nu-\delta_0$.  Consider the 
representation $T$ of $\R$ given by translation by $t$.  
Then:\\
$\mu$ is weakly measurable; \\
$\|\mu\|>0$;\\
$T_t(\mu)=T_t(\nu-\delta_0)=\nu-\delta_t$;\\
for every $A\in \Sigma$,  $T_t(\mu)(A)=0$ 
for almost all $t\in \R$.\\
It now follows from Proposition \ref{weakhypa} that the 
representation $T$ is not sup path attaining.\\
(b)  Let $\alpha$ be a real number and let 
$\Sigma,\ \mu, \nu, \delta_t$, and $T_t$ 
have the same meanings as in (a).  Define a 
representation $T^\alpha$ by
$$T_t^\alpha = e^{i\alpha t}T_t.$$ 
Arguing as in (a), it is easy to see that $T^\alpha$ 
is not sup path attaining.
}
\label{exnohypa}
\end{ex}
%%%%%%%%%%%%%%%%%%%%
%%%%%%%%%%%
%%%%%%%%%%%%
%%%%%%%%%%%
%%%%%%%%%%%%%%%%%%%%%%%%%%%%%%%%%%%%%%%%%%%%%%%%%
%%%%%%%%%%%%%%%%%%%%%%%%%%%%%%%%%%%%%%%%%%%%%%%%%%%%
%%%%%%%%%%%%%%%%%%%%%%%%%%%%%%%%%%%%%%%%%%%%%%%%%%%%%%%%%%%%%%%%%%%%%%%%%%%%%%%
%%%%%%%%%%%%%%%%%%%%%%%%%%%%%%%%%%%%%%%%%%%%%%%%%%%%%%%%%%%%%%%%%%%%%%%%%%%%%%%%

\section{Proof of the Main Theorem}

%%%%%%%%%%%
%%%%%%%%%%%%
%%%%%%%%%%%
%%%%%%%%%%%%
Throughout this section, $G$ denotes a locally compact abelian group with dual group $\Gamma$; $M(\Sigma)$ is a space of measures on a set
$\O$; and $T=(T_t)_{t\in G}$ is a sup path attaining representation of $G$ by isomorphisms of $M(\Sigma)$.   

If 
$\phi$ is in $\Linfg$, write $\left[ \phi\right]$ 
for the smallest weak-* closed translation-invariant subspace of $\Linfg$ 
containing $\phi$, and let $\cI([\phi])=\cI (\phi)$ 
denote the closed translation-invariant ideal in $L^1(G)$:\\
$$ \cI (\phi)=\{f\in L^1(G): f*\phi=0\}.$$
%%%%%%%%%%%%%%%%%%%%%%%%%%%%%%%%%%%%%%%%%%%%%%%%%%%%%%%%
It is clear that 
$ \cI (\phi)=\{f\in L^1(G): f*g=0, \forall g\in \left[\phi\right]\}$. 
%%%%%
%%%%%
%%%%%%%%
The spectrum of $\phi$, denoted by 
$\sigma \left[\phi\right]$, 
is the set of all continuous characters of $G$ 
that belong to $\left[\phi\right]$.  
This closed subset of $\Gamma$ is also given by 
%%%%%%%%%%
%%%%%%%%%%
%%%%%%%%%%
\begin{equation}
\sigma \left[\phi\right]=Z(\cI(\phi)). 
\label{spec}
\end{equation}
%%%%%%%%%%
%%%%%%%%%%
%%%%%%%%%%
(See \cite[Chapter 7, Theorem 7.8.2, (b), (c), and (d)]{rudin}.)  
%%%%%%%%%%%%%%%%
\begin{rem}
{\rm  
(a)  For a weakly measurable $\mu\in M(\Sigma)$, since $\cI(\mu)$ is a closed ideal
in $L^1(G)$, it is translation-invariant, by
\cite[Theorem 7.1.2]{rudin}.
It follows readily that 
for all $t\in G$,
%%%%%%%%%%%%%%%%%%%%
$$\cI(T_t\mu)=\cI(\mu),$$
and hence
\begin{equation}
{\rm spec}_T (T_t(\mu))={\rm spec}_T (\mu).
\label{feb.4.95}
\end{equation}
(b)  Let $\mu\in M(\Sigma)$ be weakly measurable and let $E\in \Sigma$. 
It is clear that $\cI(\mu)\supset\cI(t\mapsto (T_t\mu)(E))$.  Thus
$\sigma\left[ t\mapsto (T_t\mu)(E)\right]\subset  \spec_T(\mu)$.

}
%%%%%%%%%%%%%%%%%%%%%%%%%%%%%%%%%%%%%%%%%
%%%%%%%%%%%%%%%%%%%%%%%%%%%%%%%%%%%%%%%%%%%%%%%%%
\label{rem1/28.6}
\end{rem}
%%%%%%%%%%%%%%%%%%%%%%%%%%%%%%%
%%%%%%%%%%%%%%%%%%%%%%%%%%%%%%%%%%
%%%%%%%%%%%%%%%%%%%%%%%%%%%%%%%%%%
%%%%%%%%%%%%%%%%%%%%%%%%%%%%%%%%%%

\begin{lemma}
Let $\mu\in \cM_T(\Sigma)$.\\
(a)  If $g\in L^1(G)$ and $E\in \Sigma$, then
$$\int_G g(t)T_t\mu(E)\overline{\chi(t)}\,dt=0,$$
for all $\chi$ not in ${\rm supp}(\widehat{g})+\spec_T(\mu)$.\\
(b)\  If $\nu\in M(G)$, then 
$$\spec_T(\nu*_T\mu)\subset {\rm supp}(\widehat{\nu})\cap   \spec_T(\mu)\,.$$
(c)\  If $(k_\alpha)$ is an approximate identity for $L^1(G)$, then 
$$\|\mu\|_{M(\Sigma)}/C\leq \liminf_\alpha \left\|
t\mapsto \|k_\alpha *_T T_t\mu\|_{M(\Sigma)}\right \|_{L^\infty(G)}, $$
where $C$ is as in (\ref{ineqhypa}).
\label{ref1}
\end{lemma}
{\bf Proof.}\ (a)\ Fix $E\in\Sigma$ and $\chi\in\Gamma$, and let
$f(t)=(T_{-t}\mu)(E)\chi(t)$.
Then $\sigma\left[ f\right]\subset \chi -\spec_T(\mu)$.
Hence, for $g\in L^1(G)$ with 
$\mbox{supp}(\widehat{g})\cap (\chi-\spec_T\mu)=\emptyset$,
we have
$$\int_Gg(t)(T_t\mu)\overline{\chi}(t)dt=(g*f)(0)=0\,.$$
(b)  Immediate from $\cI(\nu*_T\mu)\supset \cI(\nu)\cup\cI(\mu)\,.$

\noindent
(c)  Fix $h$ on $\Omega$ with $|h|\leq 1$ and let $f(t)=\int_\Omega
h(x)d(T_t\mu(x))$.  Then $f\in L^\infty(G)$ and so, if $(k_\alpha)$ is an approximate
identity for $L^1(G)$, then 
$k_\alpha *f\rightarrow f$ in the weak $*$-topology of $L^\infty(G)$.  Hence
\begin{eqnarray*}
\|f\|_\infty\leq \liminf\|k_\alpha *f\|_\infty&=&
\liminf\left\|t\rightarrow \int_\Omega h(x)d(k_\alpha*_T T_t\mu) \right\|_\infty\\
&\leq &
\liminf\left\|t\rightarrow \|k_\alpha*_T T_t\mu\| \right\|_\infty\,.
\end{eqnarray*}
The proof is completed by taking the sup of $\|f\|_\infty$ 
over all $h$ and using (\ref{ineqhypa}).

\begin{lemma}
Let $f:\ G\rightarrow M(\Sigma)$ be bounded and continuous, and let $\nu\in M(G)$.  
Suppose that\\
(i)\ for all $E\in\Sigma$, $t\mapsto f(t)E$ is in $L^1_S(G)$;\\
(ii)\ for all $g\in L^1_S(G),\ \|\nu*g\|_1\leq \|g\|_1\,.$\\
Then 
\begin{equation}
 \int_G\|(\nu*f)(t)\|_{M(\Sigma)}dt\leq \int_G
\|f (t)\|_{M(\Sigma)}dt\,.
\label{ref-eq-1}
\end{equation} 
\label{ref2}
\end{lemma}
{\bf Proof.}  We suppose throughout the proof that 
the right side of (\ref{ref-eq-1}) is finite; otherwise there is nothing to prove.
Write $f_E(t)=f(t)E$.  Thus
$$(\nu*f)(t)E=(\nu*f_E)(t)$$
is continuous on $G$.  By (i) and (ii),
\begin{equation}
\int_G|(\nu*f)(t)E|dt\leq \int_G|f(t)E|dt,
\label{x1}
\end{equation}
for all $E\in\Sigma$.  Now let $P=\{E_j\}$ be a finite measurable partition of 
$\Omega$, and define 
$$h(P,t)= \sum_j|(\nu*f)(t)E_j|\,.$$
Then $h(P,\cdot)$ is continuous on $G$ and, by (\ref{x1}),
\begin{equation}
\int_G h(P,t)dt=\sum_j\int_G|(\nu*f)(t)E_j|dt
\leq \int_G\|f(t)\|_{M(\Sigma)}dt\,.
\label{x2}
\end{equation}
If $R$ is a common refinement of $P$ and $Q$, then, for all $t$,
$$\max\{h(P,t),h(Q,t)\}\leq h(R,t)\,.$$
It follows from \cite[(11.13)]{hr1} that
\begin{eqnarray*}
\int_G\|\nu*f\|_{M(\Sigma)}dt&=&\int_G \sup_Ph(P,t)dt\\
&=&\sup_P\int_G h(P,t)dt \leq \int_G\|f(t)\|_{M(\Sigma)}dt
\end{eqnarray*}
by (\ref{x2}). \\

\noindent
{\bf Proof of Theorem \ref{maintransferenceth}.} 
Let $f(t)=T_t\mu$ and suppose first that 
$f$ is continuous and that there is a neighborhood $V$
of $0\in \Gamma$ such that $V+\spec_T(\mu)\subset S$.
Choose a continuous $g\in L^1(G)$
such that
$$g\geq 0;$$
$$\widehat{g}(0)= \int_G g(t)dt=1;$$
and
$${\rm supp}(\widehat{g})\subset V\,. $$
Then, for all $E\in\Sigma$, the mapping 
$t\mapsto g(t)f(t)E$ is bounded and continuous on $G$, and, by Lemma \ref{ref1},
belongs to
$$L^1_{{\rm supp}(\widehat{g})+\spec_T(\mu)}\subset L^1_{V+\spec_T(\mu)}\subset L^1_S\,.$$ 
Hence Lemma \ref{ref2} implies that
\begin{equation}
\|\nu*(gf)\|_{L^1(G,M(\Sigma))}\leq \|gf\|_{L^1(G,M(\Sigma))}\,.
\label{star}
\end{equation}
We now proceed to show that (\ref{conctransference2}) is 
a consequence of (\ref{star}).  
We start with the right side of (\ref{star}):
%%%%%%%%%%%
%%%%%%%%%%%%%%%%%%%%%%%
%%%%%%%%%%%%
%%%%%%%%%%%
\begin{eqnarray}
\|gf\|_1         &=&     \int_Gg(t)\|f(t)\|dt \nonumber\\
		&\leq& c \|\mu\| \int_Gg(t) dt =c \|\mu\|.
\label{ohwell}
\end{eqnarray}
%%%%%%%%%%%
%%%%%%%%%%%%%%%%%%%%%%%
%%%%%%%%%%%%
%%%%%%%%%%%
Consider now the left side of (\ref{star}).  We have
%%%%%%%%%%%%%%%%%%%%%%%
%%%%%%%%%%%%
%%%%%%%%%%%
\begin{eqnarray}
\|\nu *(gf)\|_1            &=&     \int_G\|\nu *(gf)(t)\|dt \nonumber\\
		&=&     \int_G\|\int_G (gf) (t-s)d\nu(s)\| dt \nonumber\\
		&=& \int_G\|\int_G g (t-s) T_{t-s} \mu 
				d \nu (s)\| dt \nonumber\\
		&\geq&  \frac{1}{c} 
		\int_G\|\int_G T_{-t}\left[ g (t-s) T_{t-s} \mu 
				  \right]
				d \nu (s)\| dt \nonumber\\
		&=&     \frac{1}{c} 
		\int_G\|\int_G  g (t-s) T_{-s} \mu 
				d \nu (s)\| dt \nonumber\\
		&\geq&  \frac{1}{c} 
		\| \int_G \int_G  g (t-s) d t 
		 T_{-s} \mu 
				d \nu (s)\| \nonumber\\
		&=&     \frac{1}{c} \widehat{g}(0)
		\| \int_G  
		 T_{-s} \mu 
				d \nu (s)\| \nonumber\\
		&=&     \frac{1}{c}
		\| \nu*_T \mu \|.           
\label{ohwellwell}
\end{eqnarray}
%%%%%%%%%%%
%%%%%%%%%%%%%%%%%%%%%%%
%%%%%%%%%%%%
%%%%%%%%%%%
Inequalities (\ref{star}), (\ref{ohwell}), and (\ref{ohwellwell}) imply that
\begin{equation}
\|\nu*_T\mu\|\leq c^2 \|\mu\|\,.
\label{dag2}
\end{equation}
Now fix $\epsilon>0$ and let $k\in L^1(G)$ be such that
$\widehat{k}$ has compact support in $\Gamma$.  
We have by Lemma \ref{ref1}
$$\spec_T (k*_T \mu)\subset \mbox{supp}(\widehat{k}) \cap \spec_T(\mu)\,.$$
Moreover, the set 
$\{\gamma\in\Gamma:\ \int_G|1-\gamma|d|\nu|<\epsilon\}$
is a neighborhood of $0\in\Gamma$.  Hence, by the hypothesis on $S$,
there is a neighborhood $V$ of 0 in $\Gamma$  and $\gamma\in\Gamma$ such that
$$\int_G|1-\gamma|d|\nu|<\epsilon\quad 
\mbox{and}\quad V+\spec_T(k*_T\mu)\subset S-\gamma\,.$$
If $h\in L^1_{S-\gamma}$, then $\gamma h\in L^1_S$; hence
$$\|(\overline{\gamma}\nu)*h\|_1=\|\nu*(\gamma h)\|_1\leq \|\gamma h\|_1=\|h\|_1, $$
for all $h\in L^1_{S-\gamma}$.
Since translation in $L^1(G)$ is continuous, it follows that 
the map $t\mapsto T_t(k*_T\mu)$ is continuous, and hence we obtain
from (\ref{dag2}),
$$\|(\overline{\gamma}\nu)*_T(k*_T\mu)\|\leq c^2\|k*_T\mu\|\,.$$
As $\|\overline{\gamma}\nu-\nu\|<\epsilon$ and $\epsilon$ is arbitrary, we 
get
$$\|k*_T(\nu*_T\mu)\|=\|\nu*_T(k*_T\mu)\|\leq c^2\|k*_T\mu\|\leq c^3\|\mu\|\,.$$
Letting $k$ run through an approximate identity of $L^1(G)$ and using 
Lemma \ref{ref1}, we obtain $\| \nu*_T\mu \|\leq Cc^3\|\mu\|$.

%%%%%%%%%%%%%%%%%%%%%%%%%%%%%%%%%%%%%%%%%%%%%%%%%%%%%%%%%%%%%%%%%%%%%%%%%%%%%%%%%
%%%%%%%%%%%%%%%%%%%%%%%%%%%%%%%%%%%%%%%%%%%%%%%%%%%%%%%%%%%%%%%%%%%%%%%%%%%%%%%%%
%%%%%%%%%%%%%%%%%%%%%%%%%%%%%%%%%%%%%%%%%%%%%%%%%%%%%%%%%%%%%%%%%%%%%%%%%%%%%%%%%
%%%%%%%%%%%%%%%%%%%%%%%%%%%%%%%%%%%%%%%%%%%%%%%%%%%%%%%%%%%%%%%%%%%%%%%%%%%%%%%%%%
%%%%%%%%%%%%%%%%%%%%%%%%%%%%%%%%%%%%%%%%%%%%%%%%%%%%%%%%%%%%%%%%%%%%%%%%%%%%%%%%%%
%%%%%%%%%%%%%%%%%%%%%%%%%%%%%%%%%%%%%%%%%%%%%%%%%%%%%%%%%%%%%%%%%%%%%%%%%%%%%%%%%%%
%%%%%%%%%%%%%%%%%%%%%%%%%%%%%%%%%%%%%%%%%%%%%%%%%%%%%%%%%%%%%%%%%%%%%%%%%%%%%%%%%%%
%%%%%%%%%%%%%%%%%%%%%%%%%%%%%%%%%%%%%%%%%%%%%%%%%%%%%%%%%%%%%%%%%%%%%%%%%%%%%%%%%%%%%%%%%%%
\section{Transference in Spaces of Analytic measures}
%%%%%%%%%%%%%%%%%%%%%%%%%%%%%%%%%%%%%
%%%%%%%%%%%%%%%%%%%%%%%%%%%%%%%%%%%%%%%%%%%
%%%%%%%%%%%%%%%%%%%%%%%%%%%%%%%%%%%%%
%%%%%%%%%%%%%%%%%%%%%%%%%%%%%%%%%%%%%%%%%%%
%%%%%%%%%%%%%%%%%%%%%%%%%%%%%%%%%%%%%
%%%%%%%%%%%%%%%%%%%%%%%%%%%%%%%%%%%%%%%%%%%%%%%%%%%%%
%%%%%%%%%%%%%%%%%%%%%%%%%%%%%%%%%

In the remainder of this paper,
we transfer a result from Littlewood-Paley theory in $H^1(\R)$
to spaces of measures on which $\R$ is acting.
We then show 
how this transferred result implies 
with ease several of the main results of Bochner \cite{bs},
de Leeuw and Glicksberg \cite{deleeuwglicksberg}, and Forelli \cite{forelli}.

We start by setting our notation.  Let
$$H^\infty(\R)=\left\{f\in L^\infty(\R):\ \int_\R f(t)g(t)\, dt=0\ \mbox{for all}\ g\in H^1(\R)\right\}\, .$$ 
Let $M(\Sigma)=M(\O,\Sigma)$ denote a space of measures and  
let $T=(T_t)_{t\in\R}$ denote a sup path attaining representation of 
$\R$ by isomorphisms of $M(\Sigma)$.
  According to Forelli \cite{forelli},
a measure $\mu\in M(\Sigma)$ 
is called $T$-analytic if $\spec_T(\mu)\subset [0,\infty)$.
We now introduce another equivalent definition.
%%%%%%%%%%
\begin{defin}
Suppose that $T$ is a sup path attaining representation of $\R$
by isomorphisms of $M(\Sigma)$.  A measure $\mu\in \cMT$ is 
called weakly analytic if the mapping $t\mapsto T_t\mu(A)$ 
is in $H^\infty(\R)$ for every $A\in\Sigma$.
\label{weakanalyticmeasure}
\end{defin}
%%%%%%%%%%%%%%%%%%%%
It is easy to see that if $\mu$ is $T$-analytic then it is 
weakly analytic.  The converse is also true.  The proof is based on the 
fact that $[0,\infty)$ is a set of spectral synthesis (see \cite[Proposition 1.7]{ams2}).

%%%%%%%%%%%%%%%%%%%%%%%%%%%%%%%%%%%%%%%%%%%%%%%%%%%%%%%
%%%%%%%%%%%%%%%%%%%%%%%%%%%%%%%%%%%%%%%%%%%%%%%%%%%%
 For $n\in\Z$,
let $m_n\in L^1(\R)$ be the function whose Fourier
transform is piecewise linear and satisfies
%%%%%%%%%%%%%
%%%%%%%%%%%%%
%%%%%%%%%%%%%%%%
\begin{equation}
\widehat{m_n}(s)=\left\{
\begin{array}{ll}
0 & \mbox{if $s\not\in [2^{n-1},2^{n+1}]$;}\\
1 &  \mbox{if $s=2^n$.}
\end{array}
\right.
\label{mnHat}
\end{equation}
%%%%%%%%%%%%%
%%%%%%%%%%%%%
%%%%%%%%%%%%%%%%
Let $h\in L^1(\R)$ be the function whose Fourier transform is piecewise
linear and satisfies
%%%%%%%%%%%%%%%%%%%%
\begin{equation}
\widehat{h}(s)=\left\{
\begin{array}{ll}
0 & \mbox{if $s\not\in [-1,1]$,}\\
1 &  \mbox{if $|s|\leq \frac{1}{2}$.}
\end{array}
\right.
\label{anotherhhat}
\end{equation}
%%%%%%%%%%%%%
%%%%%%%%%%%%%
%%%%%%%%%%%%%%%%
It is easy to check that
%%%%%%%%%%%%%
%%%%%%%%%%%%%
\begin{equation}
\widehat{h}(s)+\sum^\infty_{n=0} \widehat{m_n}(s)
=\left\{
\begin{array}{ll}
1 & \mbox{if $s\geq -\frac{1}{2}$,}\\
0 &  \mbox{if $s\leq -1 $,}
\end{array}
\right.
\label{newhat}
\end{equation}
%%%%%%%%%%%%%
%%%%%%%%%%%%%
and that the left side of 
(\ref{newhat}) is continuous and piecewise linear.
The following theorem is a consequence of 
standard facts from Littlewood-Paley theory.
We postpone its proof to the end of this section.
%%%%%%%%%%%%%%%%%%%%%%%%
%%%%%%%%%%%%%%%%%%%%%%%%
%%%%%%%%%%%%%%%%%%%%%%%%
%%%%%%%%%%%%%%%%%%%%%%%%
%%%%%%%%%%%%%%%%%%%%%%%%

\begin{thm}
(i)  Let $h$ and $m_n$ be as above, 
$f$ be any function in $H^1(\R)$,
 and $N$ be any
nonnegative integer.  Then 
there is a 
positive constant $a$, independent of $f$ and $N$, such that
\begin{equation}
     \|h*f+ \sum^N_{n=0}
\epsilon_n m_n*f\|_1 \leq a \|f\|_1
\label{h1-conv}
\end{equation}
for any choice of $\epsilon_n=-1$ or $1$.
  (ii)  For
$f\in H^1(\R)$, we have
$$h*f+ \lim_{N\rightarrow\infty} \sum^N_{n=0}
m_n*f=f$$
unconditionally in $L^1(\R)$.  \\
%%%%%%%%%%%%%%%%%%%
%%%%%%%%%%%%%%%%%%%%%%%%
%%%%%%%%%%%%%%%%%%%%%%%%
(iii) For $f\in H^\infty (\R)$, we have
\begin{equation}
h*f+ \lim_{N\rightarrow\infty} \sum^N_{n=0}
m_n*f=f
\label{hinfty-conv}
\end{equation}
almost everywhere on $\R$.
\label{conv-hinfty}
\end{thm}

Our main theorem is the following.
%%%%%%%%%%%%%%%%%%%%%%%%%%%%%%%%%%%%%%%%%%%%%%%%%%%%%
%%%%%%%%%%%%%%%%%%%%%%%%%%%%%%%%%%%%%%%%%%%%%%%%%%%%%
%%%%%%%%%%%%%%%%%%%%%%%%%%%%%%%%%%%%%%%%%%%%%%%%%%%%%%
%%%%%%
%%%%%%%%%%%%%%%%%%%%%%%%%%%%%%%%%%%%%%%%%%%%%%%%%%%%%%%
%%%%%%%%%%%%%%%%%%%%%%%%%%%%%%%%%%%%%%%%%%%%%%%%%%%%%
%%%%%%%%%%%%%%%%%%%%%%%%%%%%%%%%%%%%%%%%%%%%%%%%%%%%%
%%%%%%%%%%%%%%%%%%%%%%%%%%%%%%%%%%%%%%%%%%%%%%%%%%%%%%
%%%%%%
%%%%%%%%%%%%%%%%%%%%%%%%%%%%%%%%%%%%%%%%%%%%%%%%%%%%%%%
\begin{thm}
Let $T$ be a representation of $\R$ in $M(\Sigma)$
that is sup path attaining,
and let $h$ and $m_n, n=0,1 ,2 , \ldots$
be as in Theorem \ref{conv-hinfty}.  Suppose that
 $\mu\in M(\Sigma)$ is weakly analytic.  Then 
\begin{equation}
\|h*_T\mu + \sum_{n=0}^N\epsilon_n m_n*_T \mu\|\leq a c^3 C\|\mu\|
\label{eq-15.2.95}
\end{equation}
for any choice of $\epsilon_n\in \{-1,1\}$, 
where $a$ is as in  (\ref{h1-conv}), 
$c$ as in (\ref{uniformlybded}) and $C$ as in (\ref{ineqhypa}).
Moreover,
%%%%%%%%%%%%%%%%%%%%%%%%%%%%%%%%%%%%%%%%%%%%%%%%%%%%%%%
\begin{equation}
\mu=h*_T\mu+\sum_{n=0}^\infty m_n*_T \mu,
\label{greatdecomp}
\end{equation}
%%%%%%%%%%%%%%%%%%%%%%%%%%%%%%%%%%%%%%%%%%%%%%%%%%%%%%%
where the series converges unconditionally in 
$M(\Sigma)$.
\label{ucc-analytic}
\end{thm}
%%%%%%%%%%%%%%%%%%%%%%%%%%%%%%%%%%%%%%%%%%%%%%%%%%%%%%%
%%%%%%%%%%%%%%%%%%%%%%%%%%%%%%%%%%%%%%%%%%%%%%%%%%%%%%%
%%%%%%%%%%%%%%%%%%%%%%%%%%%%%%%%%%%%%%%%%%%%%%%%%%%%
%%%%%%%%%%%%%

\noindent
{\bf Proof.}  To prove (\ref{eq-15.2.95}), 
combine Theorems \ref{maintransferenceth} and \ref{conv-hinfty}.
Inequality (\ref{eq-15.2.95}) states that the partial sums of the series 
$h*_T \mu+\sum_{n=0}^\infty m_n*_T \mu$
are unconditionally bounded.  
To prove that they converge unconditionally, we recall
the Bessaga-Pe\l czy\'nski Theorem from \cite{bp}.  This theorem
tells us that for any Banach
space, every unconditionally bounded series is 
unconditionally convergent if and
only if the Banach space does not contain an isomorphic copy of $c_0$.
Now since
$M(\Sigma)$ is weakly complete and $c_0$ is not 
(see \cite[Chap IV.9, Theorem 3, and IV.13.9]{dunford}), we conclude that
 $M(\Sigma)$ does not contain $c_0$.
Applying the Bessaga-Pe\l czy\'nski Theorem, we infer that there is a measure $\eta\in M(\Sigma)$
such that 
%%%%%%%%%%%%%%%%%%%%%%%%%%%%%%%%%%%%%%%%%%%%%%%%%%%%
\begin{equation}
\eta=h*_T\mu+\sum_{n=0}^\infty m_n*_T \mu
\label{infer1}
\end{equation}
%%%%%%%%%%%%%%%%%%%%%%%%%%%%%%%%%%%%%%%%%%%%%%%%%%%%
unconditionally in $M(\Sigma)$.  Moreover, $\eta$ is weakly measurable,
because of (\ref{infer1}).
It remains to show that $\eta=\mu$.
By Proposition \ref{weakhypa},
it is enough to show that
for every $A\in \Sigma$, we have
%%%%%%%%%%%%%%%%%%%%%%%%%%%%%%%%%%%%%%
\begin{equation}
T_t\mu(A)=T_t\eta(A)
\label{toshow1}
\end{equation}
%%%%%%%%%%%%%%%%%%%%%%%%%%%%%%%%%%%%%%%
for almost every $t\in \R$.  
Since $\mu$ is weakly analytic, the function $t\mapsto T_t\mu(A)$
is in $H^\infty (\R)$.
By Theorem \ref{conv-hinfty} (iii),
we have
%%%%%%%%%%%%%%%%%%%%%%%%%%%%%%%%%%%%%%%%%%%%%%%%%%%%
\begin{equation}
T_t\mu(A)=h*T_t\mu(A)+\sum_{n=0}^\infty m_n* T_t \mu(A)
\label{wehave1}
\end{equation}
%%%%%%%%%%%%%%%%%%%%%%%%%%%%%%%%%%%%%%%%%%%%%%%%%%%%
for almost every $t\in \R$.  On the other hand, from
the unconditional convergence of the series in 
(\ref{infer1}), (\ref{uniformlybded}), and (\ref{commut}), it 
follows that
%%%%%%%%%%%%%%%%%%%%%%%%%%%%%%%%%%%%%%%%%%%%%%%%%%%%
\begin{equation}
T_t\eta(A)=h*T_t\mu(A)+\sum_{n=0}^\infty m_n* T_t \mu(A)
\label{wehave2}
\end{equation}
%%%%%%%%%%%%%%%%%%%%%%%%%%%%%%%%%%%%%%%%%%%%%%%%%%%%
for all $t\in \R$.  Comparing (\ref{wehave1}) and 
(\ref{wehave2}), we see that (\ref{toshow1})
holds, completing the proof of the theorem.

%%%%%%%%%%%%%%%%%%
\begin{rem}
%%%%%%%%%%%%%%%%%%%%
{\rm It is interesting to note that
 Theorem \ref{ucc-analytic}
implies the classical
 F.\ and M.\ Riesz theorem for measures defined on the real line.
To see this, consider the representation
of $\R$ acting by translation on the Banach space
of complex regular 
Borel measures on $\R$. 
It is easy to see that a regular Borel measure is analytic
if and only if its Fourier-Stieltjes transform is supported in
$[0,\infty)$. 
In this case, each term in
(\ref{greatdecomp}) belongs to $L^1(\R)$,
being the convolution of an $L^1(\R)$
function with a regular Borel measure.
Thus, the unconditional convergence of the
series in (\ref{greatdecomp}) implies that the
measure $\mu$ is absolutely continuous.
}
\label{classicalf&m}
\end{rem}

This argument provides a 
new proof of the F.\ and M.\ Riesz Theorem,
based on Littlewood-Paley theory and the result of 
Bessaga-Pe\l czy\'nski \cite{bp}.  Also, it can be used to 
prove the following version of Bochner's generalization of the F. and M. Riesz Theorem.

\begin{thm}
Suppose that $G$ is a locally compact abelian group with
dual group $\Gamma$, and $\psi:\ \Gamma \rightarrow \R$ is a 
continuous homomorphism.  Suppose that $\nu\in M(G)$ is such that,
for every real number $s$, 
$\psi^{-1}((-\infty,s])\cap \supp (\widehat{\nu})$ is compact.  
Then $\nu$ is absolutely continuous with 
respect to Haar measure on $G$.  That is, $\nu\in L^1(G)$.
\label{bochner-thm}
\end{thm}
 {\bf Proof.}\ 
Let $\phi:\ \R\rightarrow G$ denote the continuous adjoint homomorphism
of $\psi$.  Define a representation $T=(T_t)_{t\in\R}$ of $\R$ on
the regular Borel measures $M(G)$ by
$$T_t(\mu)(A)=\mu(A+\phi(t)),$$
for all $\mu\in M(G)$ and all Borel subsets $A\subset G$.
By Example \ref{exhypa2}, $T$ is sup path attaining,
and every measure $\mu\in M(G)$ is weakly measurable.
Moreover, $\mu\in M(G)$ is weakly analytic (equivalently, $T$-analytic)
if and only if $\supp (\widehat{\mu})\subset \psi^{-1}([0,\infty))$ (see
\cite{deleeuwglicksberg}).

To prove the theorem, we can, without loss of generality,
suppose that $\supp (\widehat{\nu})\subset \psi^{-1}([0,\infty))$.
Otherwise, we consider the measure $\overline{\chi} \nu$, where 
$\psi(\chi) +\psi(\supp(\widehat{\nu}))\subset [0,\infty)$.

Let $S=\psi^{-1}([0,\infty))$.  Then $S$ is a $\cT$-set.  Applying Theorem \ref{ucc-analytic}, we see that
\begin{equation}
\nu=h*_T\nu +\sum_{n=1}^\infty m_n*_T \nu,
\label{ucc}
\end{equation}
unconditionally in $M(G)$.
For $f\in L^1(\R)$, a 
straightforward calculation shows that
$\widehat{f*_T\nu}(\chi)=\widehat{f}(\psi (\chi))\widehat{\nu}(\chi)$ for all
$\chi\in\Gamma$.  
Since
$\psi^{-1}((-\infty,s])\cap \supp (\widehat{\nu})$ is compact
for every $s\in \R$, it follows that
  $\supp(\widehat{h*_T\nu})$ 
and $\supp(\widehat{m_n*_T\nu})$
are compact.  Thus $h*_T\nu$ and $m_n*_T\nu$ 
are in $M(G)\cap L^2(G)$, and hence they belong to $L^1(G)$.
As a consequence, (\ref{ucc}) implies that 
$\nu\in L^1(G)$.\\

%%%%%%%%%%%%%%%%%%%%%%%%%%%%%%%%%%%%%%%%%%%%%%%%%%%%%
%%%%%%%%%%%%%%%%%%%%%%%%%%%%%%%%%%%%%%%%%%%%%%%%%%%%%%
With Theorem \ref{ucc-analytic} in hand, we can derive with ease
several fundamental properties of analytic measures
that were obtained previously by 
 de Leeuw-Glicksberg
\cite{deleeuwglicksberg}, and
Forelli \cite{forelli}.  We note however, that the techniques
in \cite{deleeuwglicksberg} and \cite{forelli} do not apply
in our more general settings.
%%%%%%
%%%%%%%%%%%%%%%%%%%%%%%%%%%%%%%%%%%%%%%%%%%%%%%%%%%%%%%
\begin{thm}
Let $T$ be a representation of $\R$ in $M(\Sigma)$
that is sup path attaining,
and let $\mu$ be a
weakly analytic measure in $ M(\Sigma)$.
Then the mapping $t\mapsto T_t \mu$ is 
continuous from $\R$ into $M(\Sigma)$.
\label{general-forelli}
\end{thm}
%%%%%%%%%%%%%%%%%%%%%%%%%%%%%%%%%%%%%%%%%%%%%%%%%%%%%%%
%%%%%%%%%%%%%%%%%%%%%%%%%%%%%%%%%%%%%%%%%%%%%%%%%%%%%%%
%%%%%%%%%%%%%%%%%%%%%%%%%%%%%%%%%%%%%%%%%%%%%%%%%%%%
{\bf Proof.}  Using the uniform continuity
of translation in $L^1(\R)$, it is a simple matter to show that
for any function $f\in L^1(\R)$, and any weakly measurable
$\mu\in M(\Sigma)$, the mapping
$t\mapsto f*_T T_t\mu$ is continuous.
Now use Theorem \ref{ucc-analytic} to complete the proof.\\

%%%%%%%%%%%%%%%%%%%%%%%%%%%%%%%%%%%%%%%%%%%%%%%%%%%%
%%%%%%%%%%%%%%%%%%%%%%%%%%%%%%%%%%%%%%%%%%%%%%%%%%%%%%
%%%%%%
Theorem \ref{general-forelli} is very specific to
representations of $\R$ or $\T$, in the sense that
no similar result holds on more general groups.
To see this, consider the 
group $G=\T\times \T$ with a lexicographic order on the dual
group $\Z\times \Z$.  Let $\mu_0$ denote the normalized Haar
measure on the subgroup $\{(x,y):\ y=0\}$, and consider the measure
$e^{-i x}\mu_0$.  Its spectrum is supported on
the coset $\{(m,1):\ m\in \Z\}$ and thus it is
analytic with respect to the
regular action of $G$ by translation in
$M(G)$.  Clearly, the measure $e^{-i x}\mu_0$
does not translate continuously, and so
a straightforward analog of 
Theorem \ref{general-forelli}
fails in this setting.

%%%%%%%%%%%%%%%%%%%%%%%%%%%%%%%%%%%%%%%%%%%%%%%%%%%%
%%%%%%%%%%%%%%%%%%%%%%%%%%%%%%%%%%%%%%%%%%%%%%%%%%%%%
The following application concerns bounded operators $\cP$
from $M(\Sigma)$ into $M(\Sigma)$
that commute with $T$ in the following sense:
$$\cP\circ T_t=T_t\circ \cP$$
for all $t\in \R$.  
%%%%%%%%%%%%%%%%%%%%%%%%%%%%%%%%%%%%%%%%%%%%%%%%%%%%%%%
%%%%%%%%%%%%%%%%%%%
%%%%%%%%%%%%%%%%%%
%%%%%%%%%%%%%%%%%%%%
\begin{thm}
Suppose that $T$ is a representation of $\R$ that is sup path
attaining,
and that $\cP$ commutes with $T$.
Let $\mu\in M(\Sigma)$ be weakly analytic.
Then $\cP \mu$ is also weakly analytic.
\label{aboutPforR}
\end{thm}
%%%%%%%%%%%%%%%%%%
%%%%%%%%%%%%%%%%%
{\bf Proof.}  
First note that by Theorem \ref{general-forelli}, the mapping
$t\mapsto T_t \mu(A)$
is continuous, and hence measurable.  

Now suppose that $g \in H^1(\R)$.  Again, by Theorem \ref{general-forelli},
the
map $t\mapsto g(t) T_t\mu$ is Bochner integrable.
Let
$$ \nu = \int_\R g(t) T_t \mu dt .$$
Then by properties of the Bochner integral, 
and since $\mu$ is weakly analytic, we have that
for all $A \in \Sigma$
$$ \nu(A) = \int_\R g(t) T_t\mu(A) dt = 0 .$$
Hence $\nu = 0$.
Therefore, for all $A \in \Sigma$ we have
$$ \int_\R g(t) T_t(\cP\mu)(A) dt
   =
   \int_\R g(t) \cP(T_t\mu)(A) dt
   =
   \cP\nu(A) = 0. $$
Since this is true for all $g \in H^1(\R)$, it follows that
$\cP\mu$ is weakly analytic.

%%%%%%%%%%%%%%%%%%%%%%%%%%%%%%%%%%%%%%%%%%%%%%%%%%%%%%%
\begin{defin}
Let $T$ be a sup path attaining 
representation of $G$ in $M(\Sigma)$.
A weakly measurable $\sigma$ in $M(\Sigma)$ is
called quasi-invariant if 
$T_t\sigma$ and $\sigma$
are mutually absolutely continuous for all $t\in G$.  Hence
if $\sigma$ is quasi-invariant
and $A\in \Sigma$, then 
$|\sigma|(A)=0$ if and only if $|T_t(\sigma)|(A)=0$
for all $t\in G$.
\label{qi}
\end{defin}

 %%%%%%%%%%%%%%%%%%
%%%%%%%%%%%%%%%%%%

We can use Theorem \ref{aboutPforR} to generalize a result of
de Leeuw-Glicksberg \cite{deleeuwglicksberg} 
and Forelli \cite{forelli}, concerning quasi-invariant measures.  In this application, it 
is necessary to restrict
to sup path attaining representations given by isometries of  $M(\Sigma)$.
We need a lemma.

\begin{lemma}
Suppose that $T$ is a linear isometry of $M(\Sigma)$ onto itself.
Let $\mu$, $\nu\in M(\Sigma)$.  Then,\\
(a) $\mu$ and $\sigma$ are mutually singular (in symbols, $\mu\perp\sigma$)
if and only if
$T\mu\perp T\sigma$;\\
(b)  $\mu<<\sigma$  if and only if
$T\mu<<T\sigma.$ 
\label{finallem}
\end{lemma}

\noindent
{\bf Proof.}\ For (a), simply recall that two measures $\mu$ and $\sigma$
are mutually singular if and only if $\|\mu+\sigma\|=\|\mu\|+\|\sigma\|$, and 
$\|\mu-\sigma\|=\|\mu\|+\|\sigma\|$.
For (b), it is clearly enough to prove the implication in one direction.
So suppose that 
$\mu<<\sigma$ and write $\mu=\mu_1+\mu_2$ where 
$T\mu_1<<T\sigma$ and 
$T\mu_2\perp T\sigma$.
Then 
$T\mu_1\perp T\mu_2$.  Hence 
$ \mu_1\perp  \mu_2$, and hence $\mu_2<<\mu<<\sigma$.
But $T\mu_2\perp T\sigma$ implies that $\mu_2\perp \sigma$.  So
$\mu_2=0$.  Thus $T\mu=T\mu_1<<T\sigma$.

%%%%%%%%%%%%%%%%%%
%%%%%%%%%%%%%%%%%%
%%%%%%%%%%%%%%%%%%
%%%%%%%%%%%%%%%%%%
\begin{thm}
Suppose that $T$ is a sup path attaining representation
of $\R$ by isometries of $M(\Sigma)$.  Suppose 
that $\mu\in M(\Sigma)$ is weakly analytic, and
$\sigma$ is quasi-invariant.  Write
$\mu=\mu_a+\mu_s$ for the Lebesgue decomposition of $\mu$
with respect to $\sigma$.  Then both
$\mu_a$ and $\mu_s$ are weakly analytic.  In particular,
the spectra of $\mu_a$ and $\mu_s$ are
contained in $[0,\infty)$.
\label{lebesgue-decomp-forR}
\end{thm}
%%%%%%%%%%%%%%%%%%%%%%%%%%%%%%%%%%%%%%%%%%%%%%%%%%%%
%%%%%%%%%%%%%%%%%%%%%%%%%%%%%%%%%%%%%%%%%%%%%%%%%%%%%
%%%%%%%%%%%%%%%%%%%%%%%%%%%%%%%%%%%%%%%%%%%%%%%%%%%%%%
%%%%%%%%%%%%%%%%%%
%%%%%%%%%%%%%%%%%
{\bf Proof.}  Let $\cP(\mu)=\mu_s$.  
Since $\sigma$ is quasi-invariant, the operator
$\cP$ commutes with $T$ by Lemma \ref{finallem}.  
Now apply Theorem 
\ref{aboutPforR}.
%%%%%%%%%%%%%%%%%%%%%%%%%%%%%%%%%
%%%%%%%%%%%%%%%%%%%%
%%%%%%%%%%%%
%%%%%%%%%%%

Let us finish with an example to show that the hypothesis of 
sup path attaining is required in these results.  The
next example is a variant of Example \ref{exnohypa}.

\begin{ex}
{\rm Let $\Sigma_1$ denote the sigma algebra of countable
and co-countable subsets of $\R$, let $\Sigma_2$ denote
the Borel subsets of $\R$, and let $\Sigma = \Sigma_1 \otimes
\Sigma_2$ denote the product sigma algebra on $\R \times \R$.
Let $\nu_1:\Sigma_1 \to \{0,1\}$ be the measure that takes countable sets
to $0$ and co-countable sets to $1$, let $\delta_t:\Sigma_1\to\{0,1\}$
be the measure that takes sets to $1$ if they contain $t$, and to $0$ otherwise,
and let $\nu_2$ denote the measure on $\Sigma_2$ given by
$$ \nu_2(A) = \int_A \exp(-x^2) \, dx .$$
Let $\nu = \nu_1\otimes \nu_2$, 
let $\theta = \delta_0\otimes\nu_2$, and
let $\mu = \nu-\theta$.
Finally, let $T_t$ be the representation given by 
$T_t(x,y) = (x+t,y+t)$.

Then, we see that $\nu$ is quasi-invariant, and that
$\theta$ and $\nu$ are mutually singular.
Arguing as in
Example \ref{exnohypa}, we see that $\mu$ is weakly analytic.
However, the singular part of $\mu$ with respect to $\nu$
is $-\theta$, and it may be readily seen that this is
not weakly analytic, for example
$$ T_t\theta(\R\times[-1,1])
   = \int_{-1}^1 \exp(-(x-t)^2) \, dx ,$$
is not in $H^\infty(\R)$.}
\end{ex}
%%%%%%%%%%%%%
%%%%%%%%%%%%%%%%%%%%%%%%%%%%%%%%%%%%%%%%%%%%%%%%%%%%%%%%%%%%%%%%%%%%%%%%%%
%%%%%%%%%%%%%%%%%%%%%%%%%%%%%%%%%%%%%%%%%%%%%%%%%%%%%%%
%%%%%%%%%%%%%%%%%%%%%%%%%%%%%%%%%%%%%%%%%%%%%%%%%%%%
%%%%%%%%%%%%%%%%%%%%%%%%%%%%%%%%%%%%%%%%%%%%%%%%%%%%%
%%%%%%%%%%%%%%%%%%%%%%%%%%%%%%%%%%%%%%%%%%%%%%%%%%%%%%
%%%%%%
%%%%%%%%%%%%%%%%%%%%%%%%%%%%%%%%%%%%%%%%%%%%%%%%%%%%%%%

\

We end this section by proving Theorem \ref{conv-hinfty}.
We have
\begin{equation}
\sum^\infty_{n=-\infty} \widehat{m_n}(s)
=\left\{
\begin{array}{ll}
1 & \mbox{if $s>0$;}\\
0 &  \mbox{if $s\leq 0$.}
\end{array}
\right.
\label{SummnHat}
\end{equation}
%%%%%%%%%%%%%%%%
Recall the Fej\' er kernels $\{k_a\}_{a>0}$, where
%%%%%%%%%%%%%
%%%%%%%%%%%%%
%%%%%%%%%%%%%%%%
\begin{equation}
\widehat{k_a}(s)
=\left\{
\begin{array}{ll}
1-\frac{|s|}{a} & \mbox{if $|s|<a$;}\\
0 &  \mbox{otherwise.}
\end{array}
\right.
\label{FourierTransformOfFejer}
\end{equation}
%%%%%%%%%%%%%
%%%%%%%%%%%%%
%%%%%%%%%%%%%%%%
By Fourier inversion, we see that
%%%%%%%%%%%%%
%%%%%%%%%%%%%
%%%%%%%%%%%%%%%%
\begin{equation}
m_n(x) =\exp(i 2^n x) k_{2^{n-1}}(x)+
\frac{1}{2}  \exp (i 3\ 2^{n-1} x) k_{2^{n-1}} (x).
\label{mnInTermsOfFejer}
\end{equation}
%%%%%%%%%%%%%
%%%%%%%%%%%%%
%%%%%%%%%%%%%%%%
%%%%%%%%%%%%%
%%%%%%%%%%%%%
%%%%%%%%%%%%%%%%
\begin{thm}
(i)  Let $f$ be any function in $H^1(\R)$, and let 
$M$ and $N$ be arbitrary positive integers.  Then 
there is a 
positive constant $a$, independent of $f, M,$ and $N$ such that
\begin{equation}
     \|\sum_{n=-M}^N\epsilon_n m_n*f\|_1\leq a \|f\|_1
\label{eq1lpth1}
\end{equation}
for any choice of $\epsilon_n=-1$ or 1.  (ii)  Moreover,
for $f\in H^1(\R)$,
$$\lim_{M,N\rightarrow\infty} \sum^N_{n=-M}
m_n*f=f$$
unconditionally in $L^1(\R)$.
\label{UnconditionalConvergenceOnR}
\end{thm}
%%%%%%%%%%%%%
%%%%%%%%%%%%%
%%%%%%%%%%%%%%%%
{\bf Proof.}  
The proof of (ii) is immediate from (i) and (\ref{SummnHat}),
by Fourier inversion.  
For part (i), 
use (\ref{mnInTermsOfFejer}), 
to write 
\begin{eqnarray*}
\sum_{n=-M}^N\epsilon_n m_n		&=&
\sum_{n=-M}^N\epsilon_n\exp(i 2^n x) k_{2^{n-1}}(x)+
\frac{1}{2}  \exp (i 3\ 2^{n-1} x) k_{2^{n-1}} (x)\\
					&=&
\sum_{n=-M}^{-1}\epsilon_n\exp(i 2^n x) k_{2^{n-1}}(x)
+\sum_{n=0}^N\epsilon_n\exp(i 2^n x) k_{2^{n-1}}(x)\\
&+&
\frac{1}{2}
\sum_{n=-M}^{-1} \epsilon_n
  \exp (i 3\ 2^{n-1} x) k_{2^{n-1}} (x) 
+ \frac{1}{2}\sum_{n=0}^N
	\epsilon_n
  \exp (i 3\ 2^{n-1} x) k_{2^{n-1}} (x)\\
					&=&
K_1(x)+K_2 (x)+K_3 (x)+K_4(x).
\end{eqnarray*}
Hence, to prove (\ref{eq1lpth1}) 
it is enough to show that there is a positive constant
$a$, independent of $f$ such that
$\|K_j*f\|_1\leq a \|f\|_1$, for $j=1,2,3,4$.
Appealing to \cite[Theorem 3, p.\ 114]{stein}, we will be done once we
establish that:
 \begin{equation}
 \left|\widehat{K_j}\right|\leq A,
 \label{condition1}
 \end{equation}
%%%%%%%%%%%%%
and the H\"{o}rmander condition
%%%%%%%%%%%%%%%%
\begin{equation}
\sup_{y>0} \int_{|x|>2y} |K_j(x-y)-K_j(x)|dx\leq B,
\label{condition3}
\end{equation}
where $A$ and $B$ are absolute constants.
Inequality (\ref{condition1}) 
holds with $A=1$, since the 
Fourier transforms of the summands defining the 
kernels $K_j$ have disjoint supports and are 
bounded by 1.  Condition (\ref{condition3}), is well-known. 
For a proof, see \cite[pp.\ 138-140, and 7.2.2, p. 142]{eg}.\\
%%%%%%%%%%%%%%%%%%%%%%%%%
%%%%%%%%%%%%%%%%%%%%%%%%%%%%%%%%%%%%%%%%%%%%%%%%%%%%%%%
%%%%%%%%%%%%%%%%%%%%%%%%%%%%%%%%%%%%%%%%%%%%%%%%%%%%%%%
%%%%%%%%%%%%%%%%%%%%%%%%%%%%%%%%%%%%%%%%%%%%%%%%%%%%%%%
%%%%%%%%%%%%%%%%%%%%%%%%%%%%%%%%%%%%%%%%%%%%%%%%%%%%
%%%%%%%%%%%%%%%%%%%%%%%%%%%%%%%%%%%%%%%%%%%%%%%%%%%%%
%%%%%%%%%%%%%%%%%%%%%%%%%%%%%%%%%%%%%%%%%%%%%%%%%%%%%
%%%%%%%%%%%%%%%%%%%%%%%%%%%%%%%%%%%%%%%%%%%%%%%%%%%%%%
  
%%%%%%%%%%%%%%%%%%%%%
%%%%%%%%%%%%%
%%%%%%%%%%%%%%%%

%%%%%%%%%%%%%%%%%%%%%%%%
%%%%%%%%%%%%%%%%%%%%%%%%
%%%%%%%%%%%%%%%%%%%%%%%%
%%%%%%%%%%%%%%%%%%%%%%%%
%%%%%%%%%%%%%%%%%%%%%%%%
%%%%%%%%%%%%%%%%%%%%%%%%%%%%%%%%%%%%%%%%%%%%%%%%
%%%%%%%%%%%%%%%%%%%%%%%%%%%%%%%%%%%%%%%%%%%%%%%%
{\bf Proof of Theorem \ref{conv-hinfty}.}  
Parts (i) and (ii) follow
as in Theorem \ref{UnconditionalConvergenceOnR},
so we only prove (iii).  
For notational convenience, let
%%%%%%%%%%%%%%%%%%%%%%%%%
\begin{equation}
\kappa_N(x)=\sum_{n=0}^N m_n(x),\quad N=0,\,1,\, 2\, \ldots
\label{kN}
\end{equation}
%%%%%%%%%%%%%%%%%%%
Let $V_{2^N}$ denote the de la Vall\'ee
Poussin kernels on $\R$ of order $2^N$.
Its Fourier transform is 
continuous, piecewise linear, 
and satisfies 
%%%%%%%%%%%%%%%%%%%%
%%%%%%%%%%%%%%%%%%%%
\begin{equation}
\widehat{V_{2^N}}(s)=\left\{
\begin{array}{ll}
0 & \mbox{if $s\not\in [-2^{n+1},2^{n+1}]$;}\\
1 &  \mbox{if $|s|\leq 2^n$.}      
\end{array}
\right.
\label{Vn-hat}
\end{equation}
%%%%%%%%%%%%%%%%%%%%%%%%
%%%%%%%%%%%%%%%%%%%%%%%%
It is well-known that $V_{n}$ is a summability kernel for 
$L^1(\R)$ and, in particular, that $V_{2^N}*f$
converges pointwise almost everywhere to $f$
for all $f\in L^p(\R)$, for $1\leq p\leq \infty$.
Thus, we will be done, if we can show that for $f\in
H^\infty(R)$,
\begin{equation}
V_{2^N}*f=h*f+ \sum^N_{n=0}  m_n*f.
\label{toshow15.feb.95}
\end{equation}
%%%%%%%%%%%%%%%%%%%%%%%%
%%%%%%%%%%%%%%%%%%%%%%%%
Write $V_{2^N}=(h+\kappa_{N})+(V_{2^N}-(h+\kappa_{N}))$,
where $\kappa_N$ is as in (\ref{kN}).
For $N\geq 1$, the Fourier transform of $(V_{2^N}-(h+\kappa_{N}))$,
vanishes on $[-\frac{1}{2},\infty)$.
Thus, for $f\in H^\infty(\R)$, we have
$(V_{2^N}-(h+\kappa_{N}))*f=0$, and so
(\ref{toshow15.feb.95}) follows, and the proof
is complete.

%%%%%%%%%%%%%%%%%%%%%%%%%%%%%%%%%%%%%%%%%%%%%%%%%%%%%%%
%%%%%%%%%%%%%%%%%%%%%%%%%%%%%%%%%%%%%%%%%%%%%%%%%%%%%%%
%%%%%%%%%%%%%%%%%%%%%%%%%%%%%%%%%%%%%%%%%%%%%%%%%%%%%%%
%%%%%%%%%%%%%%%%%%%%%%%%%%%%%%%%%%%%%%%%%%%%%%%%%

%%%%%%%%%%%%%%%%%%%%
%%%%%%%%%%%%%%%%%%
%%%%%%%%%%%%%%%%%
%%%%%%%%%%%%%%%%%%%%%%%%%%%%%%%%%%%%%%%%%%%%%%%%%%%%%
%%%%%%%%%%%%%%%%%%%%%%%%%%%%%%%%%%%%%%%%%%%%%%%%%%%%%
%%%%%%%%%%%%%%%%%%%%%%%%%%%%%%%%%%%%%%%%%%%%%%%%%%%%%%
%%%%%%
%%%%%%%%%%%%%%%%%%%%%%%%%%%%%%%%%%%%%%%%%%%%%%%%%%%%%%%
%%%%%%%%%%%%%%%%%%%%%%%%%%%%%%%%%%%%%%%%%%%%%%%%%%%%%%%
%%%%%%%%%%%%%%%%%%%%%%%%%%%%%%%%%%%%%%%%%%%%%%%%%%%%%%%
%%%%%%%%%%%%%%%%%%%%%%%%%%%%%%%%%
%%%%%%%%%%%%%%%%%%%%
%%%%%%%%%%%%
%%%%%%%%%%%
%%%%%%%%%%%%%%%%%%%%
%%%%%%%%%%%%%%%%%%
%%%%%%%%%%%%%%%%%
%%%%%%%%%%%%%%%%%%%%%%%%%%%%%%%%%%%%%%%%%%%%%%%%%%%%%
%%%%%%%%%%%%%%%%%%%%%%%%%%%%%%%%%%%%%%%%%%%%%%%%%%%%%
%%%%%%%%%%%%%%%%%%%%%%%%%%%%%%%%%%%%%%%%%%%%%%%%%%%%%%
%%%%%%

%%%%%%%%%%%%%%%%%%%%%%%%%%%%%%%%%%%%%%%%%%%%%%%%%%%%
%%%%%%%%%%%%%%%%%%%%%%%%%%%%%%%%%%%%%%%%%%%%%%%%%%%%%%
{\bf Acknowledgements}  The work of the second author was supported
by a grant from the National Science Foundation (U.S.A.).
%%%%%%%%%%%%%%%%%%%%%%%%%%%%%%%%%%%%%%%%%%%%%%%%%%%%%%%%%%%%%%%
%%%%%%%%%%%%%%%%%%%%%%%%%%%%%%%%%%%%%%%%%%%%%%%%%%%%%%%%%%%%%
%%%%%%%%%%%%%%%%%%%%%%%%%%%%%%%%%%%%%%%%%%%%%%%%%%%%%%%%%%%%%%
%%%%%%%%%%%%%%%%%%%%%%%%%%%%%%%%%%%%%%%%%%%%%%%%%%%%%%%%%%%%%%%
%%%%%%%%%%%%%%%%%%%%%%%%%%%%%%%%%%%%%%%%%%%%
%%%%%%%%%%%%%%%%%%%%%%%%%%%%%%%%%%%%%%%%%%%%%%%%%%%%%%%%%%%%%%
%%%%%%%%%%%%%%%%%%%%%%%%%%%%%%%%%%%%%%%%%%%%%%%%%%%%%%%%%%%%%%%%
%%%%%%%%%%%%%%%%%%%%%%%%%%%%%%%%%%%%%%%%%%
%%%%%%%%%%%%%%%%%%%%%%%%%%%%%%%%%%%%%%%%%%%%%%%%%%%%%%%%%%%%%%%%%%%%%%%%%%%%%%%%%%%%%%%%%%%%%%%%%%%%%%%%%%%%%%%%%%%%%%%%%%%%%%%%%%%%%%%%%%%%%%%%%%%%%%%%%%%%%%%%%%%%%%%%%%%%%%%%%%%%%%%%%%%%%%%%%%%%%%%%%%%%%%%%%%%%%%%%%%%%%%%%%%%%%%%%%%%%%%%%%%%%%%%%%%%%%%%%%%%%%%%%%%%%%%%%%%%%%%

\end{document}